\documentclass[12pt,english,sort&compress]{article}
\usepackage{mathpazo}
\usepackage[T1]{fontenc}
\usepackage[utf8]{inputenc}
\usepackage{geometry}
\usepackage{babel}
\usepackage{units}
\usepackage{textcomp}
\usepackage{mathtools}
\usepackage{enumitem}
\usepackage{amsmath}
\usepackage{amssymb}
\usepackage{graphicx}
\usepackage[numbers]{natbib}
\usepackage[unicode=true,
 bookmarks=true,bookmarksnumbered=false,bookmarksopen=false,
 breaklinks=false,pdfborder={0 0 0},pdfborderstyle={},backref=false,colorlinks=false]
 {hyperref}
\hypersetup{pdftitle={Ultra-Fast Reactive Transport Simulations When Chemical Reactions Meet Machine Learning: Chemical Equilibrium},
 pdfauthor={Allan M.M. Leal}}

\makeatletter
\newcommand{\lyxaddress}[1]{
\par {\raggedright #1
\vspace{1.4em}
\noindent\par}
}

\@ifundefined{date}{}{\date{}}
\usepackage{fouriernc}
\usepackage{textcomp} % Needed for degC
\usepackage{kbordermatrix}
\usepackage{comment}
\usepackage{graphicx}
\newcommand\wlength{2.5em}
\newcommand\w[1]{\makebox[\wlength]{$#1$}}
\newcommand\minus[1]{\mathllap{-}#1}

\usepackage{titlesec} % Allows customization of titles
\usepackage[labelfont=bf]{caption}

\makeatother

\begin{document}

\title{\vspace{-1.5\baselineskip}\textbf{Ultra-Fast Reactive Transport
Simulations When Chemical Reactions Meet Machine Learning: Chemical
Equilibrium}}

\author{Allan M. M. Leal$^{\text{a,}}$\thanks{Corresponding author}\\
{\footnotesize{}\href{mailto:allan.leal@erdw.ethz.ch}{allan.leal@erdw.ethz.ch}}
\and Dmitrii A. Kulik$^{\text{b}}$\\
{\footnotesize{}\href{mailto:dmitrii.kulik@psi.ch}{dmitrii.kulik@psi.ch}}
\and Martin O. Saar$^{\text{a}}$\\
{\footnotesize{}\href{mailto:msaar@ethz.ch}{msaar@ethz.ch}}}

\maketitle
\vspace{-1\baselineskip}

\lyxaddress{\begin{center}
{\small{}$^{\text{a}}$}\emph{\small{}Geothermal Energy and Geofluids
Group, }\\
\emph{\small{}Department of Earth Sciences, ETH Zürich, Switzerland\\[0.5em]}{\small{}$^{\text{b}}$}\emph{\small{}Laboratory
for Waste Management, Nuclear Energy and Safety Research Department,}\\
\emph{\small{}Paul Scherrer Institut, 5232 Villigen PSI, Switzerland}
\par\end{center}}

\vspace{-1\baselineskip}
\begin{abstract}
During reactive transport modeling, the computational cost associated
with chemical reaction calculations is often 10–100 times higher than
that of transport calculations. Most of these costs results from chemical
equilibrium calculations that are performed at least once in every
mesh cell and at every time step of the simulation. Calculating chemical
equilibrium is an iterative process, where each iteration is in general
so computationally expensive that even if every calculation converged
in a single iteration, the resulting speedup would not be significant.
Thus, rather than proposing a fast-converging numerical method for
solving chemical equilibrium equations, we present a machine learning
method that enables new equilibrium states to be quickly and accurately
estimated, whenever a previous equilibrium calculation with similar
input conditions has been performed. We demonstrate the use of this
smart chemical equilibrium method in a reactive transport modeling
example and show that, even at early simulation times, the majority
of all equilibrium calculations are quickly predicted and, after some
time steps, the machine-learning-accelerated chemical solver has been
fully trained to rapidly perform all subsequent equilibrium calculations,
resulting in speedups of almost two orders of magnitude. We remark
that our new \emph{on-demand machine learning method} can be applied
to any case in which a massive number of sequential\slash{}parallel
evaluations of a computationally expensive function $f$ needs to
be done, $y=f(x)$. We remark, that, in contrast to traditional machine
learning algorithms, our on-demand training approach does not require
a statistics-based training phase before the actual simulation of
interest commences. The introduced on-demand training scheme requires,
however, the first-order derivatives $\partial f/\partial x$ for
later smart predictions.

\end{abstract}

\section{Introduction}

During reactive transport simulations, the following three degrees
of chemical reactivity behavior may be observed across space and over
time: weak, moderate, and intense. For example, at relatively distant
locations from where a fluid enters a medium, none or little reactivity
may occur, potentially over long periods of time, if those regions
were initially in chemical, thermal, and mechanical equilibrium. Eventually,
these equilibrium conditions are gradually disrupted due to fluid
flow, heat transfer, and chemical transport, causing some moderate
chemical reactivity. Later, the introduced perturbations reach those
once calm locations and cause relatively intense local chemical changes.

Interestingly, after such a wave of perturbations has passed, very
often local chemical reactivity becomes weak once again. This implies
that only a relatively small percentage of the entire medium is experiencing
fast and substantial chemical changes at any moment in time. As a
result, most chemical reaction calculations are performed for low-reactivity
regions at every time step of the simulation. Even though reaction
calculations are relatively faster in low-reactivity regions, compared
to high-reactivity regions, their overall computational cost is typically
10–100 times greater than that for calculations of physical processes
(e.g., mass transport, heat transfer). Thus, reactive transport simulations
spend most of the computation time calculating chemical reactions
and not transport processes, and typically require very long compute
times. Hence, speeding up chemical calculations, without compromising
accuracy, is crucial for significant performance gains in reactive
transport simulations.

Whenever chemical reaction rates are considerably faster than rates
of physical processes, the \emph{local chemical equilibrium assumption}
is a plausible and sufficient rate model for chemical reactions. In
general, however, a combination of fast and slow reaction rates is
present, so that the fast reactions are reasonably modeled employing
the chemical equilibrium assumption and slow reaction rates are modeled
using chemical kinetics. This approach is termed the \emph{local partial
chemical equilibrium assumption} \citep{Ramshaw1980,Ramshaw1981,Ramshaw1985,Ramshaw1995,Lichtner1985,Steefel1994,Steefel1996}.
Nevertheless, what needs to be noticed about these assumptions is
that chemical equilibrium calculations are needed to model chemical
processes in either case. In fact, such equilibrium calculations need
to be performed at least once at every mesh node\slash{}cell during
every transport time step. Thus, millions to billions of chemical
equilibrium calculations tend to accumulate over the course of massive
numerical reactive transport simulations, particularly when using
fine-resolution meshes, large three-dimensional domains, and\slash{}or
long simulation times. Hence, speeding up reactive transport modeling
by a significant factor can only be accomplished by \emph{accelerating
chemical reaction calculations in general, and chemical equilibrium
calculations} in particular. Here, we focus on accelerating chemical
equilibrium calculations, whereas a companion paper (\citet{Leal2017c})
introduces a similar strategy to accelerate chemical kinetics  calculations.

Chemical equilibrium calculations are computationally expensive. This
is because \emph{(i)} they involve the iterative solution of a system
of non-linear algebraic equations, requiring at every iteration \emph{(ii)}
the evaluation of thermodynamic properties, such as activity coefficients,
concentrations, activities, chemical potentials (e.g., using the \citet{Pitzer1973}
model for aqueous solutions and the \citet{Peng1976} model for gaseous
solutions), and \emph{(iii)} the solution of a system of linear equations
with dimension on the order of either the number of chemical species
or the number of chemical elements \citep{Leal2017}. Over several
decades, major advances in developing fast, accurate, and robust methods
for chemical equilibrium calculations have been made, either based
on \emph{Gibbs energy minimization} (GEM) or \emph{law of mass action}
(LMA) formulations \citep{White1958,Smith1980,Smith1982,Alberty1992b,Crerar1975,DeCapitani1987,Eriksson1989,Ghiorso1994,Gordon1971,Gordon1994a,Gordon1994b,Harvey2013,Harvie1987,Karpov1997,Karpov2001,Karpov2002,Koukkari2011a,Kulik2013,Leal2013,Leal2014,Leal2016a,Leal2016c,Leal2017,Morel1972,Neron2012,Nordstrom1979,Trangenstein1986,VanZeggeren1970,Vonka1995,Wolery1975,Zeleznik1960}.
However, even if one could devise an algorithm that would always converge
in one single iteration, instead of the typical few to dozens of iterations,
the computational cost of chemical equilibrium calculations would
still be dominant among all other calculations. The question of utmost
importance is, thus:\emph{ Is it possible to calculate chemical equilibrium
without actually solving the non-linear equations governing equilibrium,
which requires several evaluations of thermodynamic models and as
many solutions of systems of linear equations?} In this paper, we
demonstrate that this can, indeed, be achieved with exceptional speed
and accuracy.

Our here-introduced smart chemical equilibrium algorithm operates
like a human being. A human begins life with none to little problem
solving skills. At an early age, the child is exposed to a variety
of challenges and learns how to solve them, often with assistance
from others. However, once the child has learned how to solve a specific
problem, it can solve \emph{similar problems} without requiring assistance.
Solving a problem with the help of an already acquired skill is \emph{much
faster} than having to first learn how to solve that problem, since
\emph{learning} is a difficult and time-consuming process. As the
child grows, it can solve more and more problems until, eventually,
as an adult, external help is needed only occasionally, under exceptional
circumstances, i.e., in situations that have rarely or never been
encountered before. Importantly, in contrast to going to school, this
way of ``learning by doing'' may be viewed as on-demand learning
or training, a highly efficient way of learning as one only learns
what is actually needed at least once. 

The above metaphor illustrates some key features of our new, on-demand
machine learning algorithm for fast chemical equilibrium calculations.
Initially, the algorithm has no recollection of any past calculation.
Thus, the first time the algorithm faces a new chemical equilibrium
problem, it fully solves the problem and saves the inputs and outputs
describing the calculated equilibrium state. This is, as discussed
previously, a computationally expensive process, now considered as
an act of learning. The next time the algorithm is employed, it recalls
that a calculation has previously been performed and attempts to quickly
predict the new chemical equilibrium state, using \emph{sensitivity
derivatives} to project the previously recorded equilibrium state
into the new one. Then, the algorithm checks if the prediction is
accurately acceptable within some given tolerance. If it is, then
the predicted equilibrium state is accepted as a solution, otherwise
a full chemical equilibrium calculation is performed and the corresponding
inputs and outputs recorded to describe the newly learned equilibrium
state. This way, the memory of past and fully solved equilibrium problems
grows, permitting searches among all saved and learned problems to
determine the one that is closest to the new equilibrium problem. 

In contrast to traditional machine learning algorithms, which first
require a time-consuming training phase, in which potentially many
outcomes are calculated for problems that may later never occur, our
new machine learning algorithm is only trained with problems that
actually occur, relying on the recurrence of subsequent \emph{similar
problems} to quickly respond with an accurate prediction. Hence, traditional
machine learning is equivalent to going to school, while our machine
learning algorithm is more akin to on-demand learning or ``learning
by doing.'' Furthermore, there is no clear criterion to decide, when
the training phase of traditional machine learning algorithms is ideally
terminated, i.e., when it is time to ``leave school.'' 

The here-introduced smart chemical equilibrium algorithm is particularly
useful for applications that require repetitive calculations, such
as chemical equilibrium calculations during reactive transport modeling.
For such applications, the algorithm can eventually achieve a knowledgeable
state, so that, after some time steps during a reactive transport
simulation, all equilibrium calculations can be rapidly and accurately
performed using previously learned \emph{key chemical equilibrium
states}. Even if the smart algorithm occasionally faces some exceptional
circumstances that require additional on-demand training (e.g., a
sudden change in the chemistry or temperature of the fluid entering
a region of the medium), one can still reasonably expect significant
speedups with this machine-learning-accelerated algorithm, if the
occasions during which learning is needed are considerably less frequent
than the occasions during which the algorithm can make smart predictions.

The ``intelligence'' of the here-introduced algorithm is, thus,
a combination of both the memory of already solved equilibrium problems,
that needed to be solved anyway, and its ability to ``learn'' new
ones on-demand. Physiologically speaking, the algorithm is like a
\emph{brain} that not only can record its experiences on solving chemical
equilibrium problems, but can also make decisions about when it needs
additional training, following its judgment on the accuracy\slash{}acceptability
of an estimated chemical equilibrium outcome. This on-demand training
is key to:
\begin{enumerate}[label=\emph{(\roman*)}]
\item \emph{Performance}: learning and keeping only what is needed results
in compact storage and thus fast search operations when finding the
closest past equilibrium problem already solved;
\item \emph{Accuracy}: finer control on minimizing errors resulting from
predicted equilibrium states by performing additional training whenever
needed;
\item \emph{Convenience}: neither a dedicated prior training stage nor anticipation
of possible chemical conditions during the simulation are required;
\item \emph{Simplicity}: users immediately benefit from high-performance
reactive transport simulations without having to prepare any prior
configuration and training of the smart chemical equilibrium solver. 
\end{enumerate}
Our on-demand training philosophy differs radically from most previous
efforts in speeding up chemical reaction calculations in reactive
transport simulations, which use conventional machine learning methods.
\citet{Jatnieks2016}, for example, present the steps necessary to
construct a so-called surrogate model for fast speciation calculations.
The construction of this surrogate model requires a prior training
phase, during which many random input conditions are used in a speciation
solver (PHREEQC \citep{Parkhurst2013}) and the resultant outputs
collected for statistical learning. During their numerical experiment,
32 different statistics and machine learning methods were tried to
identify the potentially best one. For a specific reactive transport
modeling problem, they collected all possible input-output combinations
in speciation calculations, and from these combinations, 7880 random
inputs were used for training the statistics model, which represented
80\% of all collected input conditions. The various constructed surrogate
models were subsequently used in a reactive transport simulation.
Each surrogate model was constructed with different 7880 input-output
samples.

Our new \emph{on-demand  machine learning  approach} has clear advantages
over conventional statistics\nobreakdash-based machine learning methods.
First, the use of sensitivity derivatives of the calculated equilibrium
states results in a machine learning method that better understands
the behavior of chemical systems regarding its reaction behavior following
changes in the input equilibrium conditions. Secondly, the use of
these sensitivity derivatives permits extrapolating and predicting
new equilibrium states with much higher accuracy and confidence. Thirdly,
it requires no statistical training before it can be applied in a
reactive transport simulation, because its on-demand training characteristics
allow it to spontaneously learn only what is needed during a reactive
transport simulation to keep the predictions accurate enough.  Furthermore,
our approach is not only simpler from a user point of view, but also
potentially faster, as it will save only a few input conditions compared
to any statistical approach that can only get more accurate the more
it knows a priori. Finally, our method produces outputs that always
satisfy the mass conservation conditions of chemical elements and
electric charge, since these constraints are incorporated into the
calculation of the sensitivity derivatives (see \citet{Leal2017})
and, thus, our method contrasts with statistics\nobreakdash-based
machine learning methods that can fail predicting equilibrium states
that satisfy given mass balance conditions.

This communication is organized as follows. In Section~\ref{sec:Definitions-and-Notation},
we introduce definitions and notation that are needed to describe
the new algorithm rigorously. In Section~\ref{sec:Method}, we formulate
the smart chemical equilibrium method, providing details on the prediction
of new equilibrium states from previously saved and learned states
as well as on the acceptance criteria for the predicted equilibrium
states. In Section~\ref{sec:Results}, we compare the performance
and accuracy of a reactive transport simulation using the here-proposed
smart equilibrium algorithm against that using a conventional Newton-based
equilibrium algorithm. Finally, we discuss in Section~\ref{sec:Conclusions-and-Future-Work}
the implications and conclusions of this study together with a planned
road-map for further research efforts in this direction. 

The here-introduced smart chemical equilibrium algorithm is implemented
in Reaktoro \citep{Leal2015b}, a unified open-source framework for
modeling chemically reactive systems (\href{http://reaktoro.org}{reaktoro.org}).

\section{Definitions and Notation\label{sec:Definitions-and-Notation}}

Consider a chemical system is a collection of \emph{chemical species}
composed of one or more \emph{components} and distributed among one
or more \emph{phases}. The species can be substances such as aqueous
ions (e.g., Na$^{+}$(aq), Cl$^{-}$(aq), HCO$_{3}^{-}$(aq)), neutral
aqueous species (e.g., SiO2(aq), CO$_{2}$(aq), H$_{2}$O(l)), gases
(e.g., CO$_{2}$(g), CH$_{4}$(g), N$_{2}$(g)), and pure condensed
phases (e.g., CaCO$_{3}$(s, calcite), SiO$_{2}$(s, quartz), Al$_{2}$Si$_{2}$O$_{5}$(OH)$_{4}$(s,
kaolinite)). The phases are each composed of one or more different
chemical species with homogeneous properties within their boundaries;
a phase can be aqueous, gaseous, liquid, solid solutions, a pure mineral,
a plasma, etc. Note that substances with the same chemical formula,
but in different phases, are distinct species (e.g., CO$_{2}$(aq)
and CO$_{2}$(g) are different species). The components are \emph{chemical
elements }(e.g., H, O, C, Na, Cl, Ca, Si) and \emph{electrical charge}
(Z), but it can also be a linear combination of these, commonly known
as \emph{primary species} (e.g., H$^{+}$(aq), H$_{2}$O(l), CO$_{2}$(aq)).
We shall use from now on the words elements and components interchangeably,
with elements not only denoting chemical elements but also electrical
charge \citep{Leal2017}.

A chemical system can exist at infinitely many \emph{chemical states}.
A chemical state is defined here as the triplet $(T,P,n)$, where
$T$ is temperature, $P$ is pressure, and ${n=(n_{1},\ldots,n_{\mathrm{N}})}$
is the vector of molar amounts of the species, with $n_{i}$ denoting
the mole amount of the $i$th species and N the number of species.
If we denote by ${b=(b_{1},\ldots,b_{\text{E}})}$ the vector of molar
amounts of the elements, with $b_{j}$ denoting the amount of the
$j$th element and $\text{E}$ the number of elements, then $b$ and
$n$ are related through the following mass conservation equation:
\begin{equation}
An=b,\label{eq:mass-balance}
\end{equation}
where $A$ is the \emph{formula matrix }of the chemical system \citep{Smith1982},
with $A_{ji}$ denoting the coefficient of the $j$th element in the
$i$th species, and thus a matrix with dimensions $\text{E}\times\text{N}$.
Below is an example of a formula matrix for a simple chemical system,
containing ${\text{N}=10}$ species and ${\text{E}=4}$ elements,
distributed among two phases, namely an aqueous and a gaseous phase,
with names of aqueous species suffixed with (aq), except H$_{2}$O(l),
and gaseous species with (g): \begin{equation}
A=\kbordermatrix{
& \mathrm{H_2O(l)} & \mathrm{H^+(aq)} & \mathrm{OH^-(aq)} & \mathrm{CO_2(aq)} & \mathrm{HCO_3^-(aq)} & \mathrm{CO_3^{2-}(aq)} & \mathrm{O_2(aq)} & \mathrm{H_2(aq)} & \mathrm{CO_2(g)} & \mathrm{H_2O(g)}\\
\mathrm{H} & \w{2} & \w{1} & \w{1} & \w{0} & \w{1} & \w{0} & \w{0} & \w{2} & \w{0} & \w{2} \\
\mathrm{O} & 1 & 0 & 1 & 2 & 3 & 3 & 2 & 0 & 2 & 1\\
\mathrm{C} & 0 & 0 & 0 & 1 & 1 & 1 & 0 & 0 & 1 & 0\\
\mathrm{Z} & 0 & 1 & \minus{1} & 0 & \minus{1} & \minus{2} & 0 & 0 & 0 & 0
}.\label{eq:formula-matrix-example}
\end{equation}In the mathematical presentation that follows, we assume, for convenience
reasons, that the formula matrix $A$, of the corresponding chemical
system, is full-rank (i.e., the rows of $A$ are linearly independent)
so that $\text{rank}(A)=\text{E}$.

\section{Method\label{sec:Method}}

Let a chemical equilibrium calculation be represented as:
\begin{equation}
n=\varphi(T,P,b).\label{eq:equilibrium-func}
\end{equation}
The chemical equilibrium function,  $\varphi$, is an abstraction
of the operations needed to solve the fundamental Gibbs energy minimization
problem:
\begin{equation}
\min_{n}G(T,P,n)\quad\text{subject to}\quad\left\{ \begin{aligned}An=b\\
n\geq0
\end{aligned}
\right.,\label{eq:gem-problem}
\end{equation}
at prescribed conditions of temperature, $T$, pressure, $P$, and
element amounts, ${b=(b_{1},\ldots,b_{\text{E}})}$. In this problem,
one seeks ${n=(n_{1},\ldots,n_{\text{N}})}$ that minimizes the total
Gibbs energy function of the system: 
\begin{equation}
G=\sum_{i=1}^{\text{N}}\mu_{i}n_{i},
\end{equation}
subject to the elemental mass conservation constraint equations, $An=b$,
and the non-negativity constraints for the species amounts, ${n_{i}\geq0}$.
The \emph{chemical potential} of the $i$th species, $\mu_{i}=\mu_{i}(T,P,n)$,
is defined as:
\begin{equation}
\mu_{i}=\mu_{i}^{\circ}+RT\ln a_{i},
\end{equation}
with $R$ denoting the universal gas constant, $\mu_{i}^{\circ}=\mu_{i}^{\circ}(T,P)$
the \emph{standard chemical potential} of the $i$th species, and
$a_{i}=a_{i}(T,P,n)$ the \emph{activity} of the $i$th species. More
details on how this chemical equilibrium problem can be solved, using
either Gibbs energy minimization (GEM) or law of mass action (LMA)
methods, are presented in Section~\ref{subsec:Chemical-equilibrium-equations}
or in more detail in \citet{Leal2016a,Leal2016c,Leal2017}.

\subsection{First-Order Taylor Approximation\label{subsec:First-order-Taylor-approximation}}

Assume a chemical equilibrium calculation has been done previously
with inputs $(T^{\text{p}},P^{\text{p}},b^{\text{p}})$, and the new
chemical equilibrium calculation needs to be performed with inputs
$(T^{\text{q}},P^{\text{q}},b^{\text{q}})$. Rather than computing
$n^{\text{q}}$ using the computationally expensive chemical equilibrium
function~$\varphi$:
\begin{equation}
n^{\text{q}}=\varphi(T^{\text{q}},P^{\text{q}},b^{\text{q}}),
\end{equation}
we want to try first to estimate $n^{\text{q}}$ using the following
\emph{first-order Taylor approximation}:
\begin{equation}
n^{\text{q}}=n^{\text{p}}+\frac{\partial\varphi^{\mathrlap{\text{p}}}}{\partial T}(T^{\text{q}}-T^{\text{p}})+\frac{\partial\varphi^{\mathrlap{\text{p}}}}{\partial P}(P^{\text{q}}-P^{\text{p}})+\frac{\partial\varphi^{\mathrlap{\text{p}}}}{\partial b}(b^{\text{q}}-b^{\text{p}}),\label{eq:smart-estimate}
\end{equation}
where $\partial\varphi^{\text{p}}/\partial T$, $\partial\varphi^{\text{p}}/\partial P$,
and $\partial\varphi^{\text{p}}/\partial b$ are \emph{sensitivity
derivatives }(with dimensions $\text{N}\times1$, $\text{N}\times1$,
and $\text{N}\times\text{E}$, respectively) of the previous chemical
equilibrium state. These sensitivity derivatives are measures of how
the species amounts change when the previous chemical equilibrium
state is perturbed by infinitesimal changes in temperature, pressure,
and the amounts of elements. They can be equivalently written as $\partial n^{\text{p}}/\partial T$,
$\partial n^{\text{p}}/\partial P$, and $\partial n^{\text{p}}/\partial b$. 

By using these sensitivity derivatives, we can quickly and accurately
estimate all chemical equilibrium states in the vicinity of some previous,
and fully known, chemical equilibrium state.

\subsection{Acceptance Testing}

\begin{figure}
\begin{centering}
\includegraphics{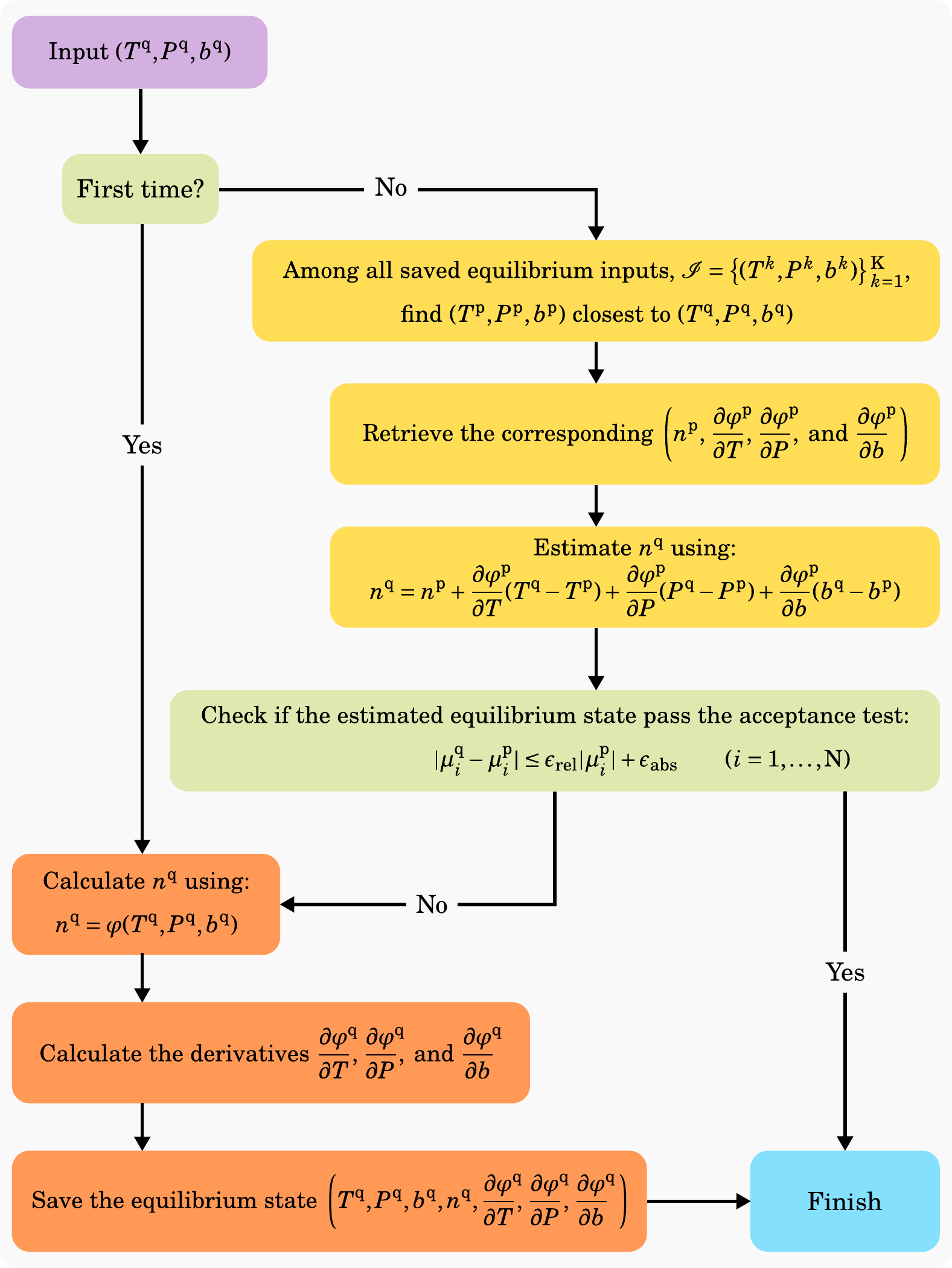}
\par\end{centering}
\caption{\label{fig:algorithm-diagram}The diagram of the proposed machine
learning algorithm for fast chemical equilibrium calculations with
a supervised, on-demand training strategy.}
\end{figure}

Once a predicted chemical equilibrium state is calculated, one needs
to test if it can be accepted. Because we are using a first-order
Taylor approximation, we need to ensure that the new estimated chemical
equilibrium state is not too far from the previous fully-calculated
equilibrium state, used as a reference point. This could be done naively
by checking how much the species amounts changed from one state to
the other using the following test condition for all species:
\begin{equation}
|n_{i}^{\text{q}}-n_{i}^{\text{p}}|\leq\epsilon_{\text{rel}}|n_{i}^{\text{p}}|+\epsilon_{\text{abs}}\qquad(i=1,\ldots,\text{N}),\label{eq:eq:acceptance-test-n}
\end{equation}
which controls how much the absolute and relative changes in the new
estimated species amounts, $n_{i}^{\text{q}}$, can be tolerated for
given absolute and relative tolerance parameters $\epsilon_{\text{abs}}$
and $\epsilon_{\text{rel}}$, respectively. 

However, there is a major disadvantage of the previous tolerance test:
it does not ``understand''  the \emph{thermodynamic behavior of
stable phases}. Consider a chemical system with two stable phases:
an aqueous solution saturated with a mineral. Clearly, adding more
of that mineral to the system will not alter the composition of the
fluid and just accumulate to a higher total amount of solids. Under
such conditions, the acceptance test based on species amounts would
fail, even though the estimated chemical equilibrium state, using
sensitivity derivatives, could be done very accurately for \emph{any
large amounts} of added mineral. By adding 1000 moles of the mineral,
and enforcing the constraint that, for example, not more than 10\%
of a species amount can vary from the previous state to the new estimated
state, that acceptance test would trigger hundreds of unnecessary
training equilibrium calculations if the initial mineral amount is
1 mol. 

For the previous example, assuming the perturbation of the system
is made only by adding the mineral, i.e., without changing temperature
or pressure, there is one thermodynamic quantity that would remain
constant: the \emph{chemical potentials of the species}, $\mu_{i}$.
Because of this behavior, we use instead the following acceptance
test in terms of chemical potentials:
\begin{equation}
|\mu_{i}^{\text{q}}-\mu_{i}^{\text{p}}|\leq\epsilon_{\text{rel}}|\mu_{i}^{\text{p}}|+\epsilon_{\text{abs}}\qquad(i=1,\ldots,\text{N}),\label{eq:acceptance-test-mu}
\end{equation}
where $\mu_{i}^{\text{q}}$ is the estimated, \emph{not evaluated},
chemical potential of the $i$th species at the new chemical equilibrium
state: 
\begin{equation}
\mu^{\text{q}}=\mu^{\text{p}}+\frac{\partial\mu^{\mathrlap{\text{p}}}}{\partial T}(T^{\text{q}}-T^{\text{p}})+\frac{\partial\mu^{\mathrlap{\text{p}}}}{\partial P}(P^{\text{q}}-P^{\text{p}})+\frac{\partial\mu^{\mathrlap{\text{p}}}}{\partial n}(n^{\text{q}}-n^{\text{p}}),\label{eq:mu-estimated}
\end{equation}
where $\partial\mu^{\text{p}}/\partial T$, $\partial\mu^{\text{p}}/\partial P$,
and $\partial\mu^{\text{p}}/\partial n$ are chemical potential derivatives
(with dimensions $\text{N}\times1$, $\text{N}\times1$, and $\text{N}\times\text{N}$,
respectively) evaluated at the previous equilibrium state. 

\textbf{Remark:} For LMA methods, in which no access to standard chemical
potentials, $\mu_{i}^{\circ}$, exists in some cases (e.g., the thermodynamic
database used only contains equilibrium constants of reactions), a
similar, but \emph{not equivalent}, test is the use of activities:
\begin{equation}
|\ln a_{i}^{\text{q}}-\ln a_{i}^{\text{p}}|\leq\epsilon_{\text{rel}}|\ln a_{i}^{\text{p}}|+\epsilon_{\text{abs}}\qquad(i=1,\ldots,\text{N}).\label{eq:acceptance-test-lna}
\end{equation}
Note, however, that activities are in general less sensitive to temperature
variations than standard chemical potentials and, thus, the above
alternative acceptance test would be less rigorous than equation~(\ref{eq:acceptance-test-mu})
and more indifferent towards temperature changes. To fix this problem,
one could use the conversion approach detailed in \citet{Leal2016b},
which permits apparent standard chemical potentials of the species
to be calculated using equilibrium constants of reactions.

\subsection{Machine Learning Chemical Equilibrium Algorithm \label{subsec:Smart-chemical-equilibrium}}

The machine learning chemical equilibrium algorithm (or, alternatively,
smart chemical equilibrium algorithm) proposed here is capable of
``remembering'' past calculations and can make use of these calculations
to quickly and accurately estimate new equilibrium states. Figure~\ref{fig:algorithm-diagram}
is a flowchart that illustrates the main steps of the algorithm.

During its first chemical equilibrium calculation, with given $(T^{\text{q}},P^{\text{q}},b^{\text{q}})$
conditions, the algorithm does so by solving the Gibbs energy minimization
problem (\ref{eq:gem-problem}). This is represented here, using the
computationally expensive equilibrium function,~$\varphi$:
\begin{equation}
n^{\text{q}}=\varphi(T^{\text{q}},P^{\text{q}},b^{\text{q}}).
\end{equation}
Once this is done, we need to save not only the input conditions $(T^{\text{q}},P^{\text{q}},b^{\text{q}})$,
but also the corresponding, calculated equilibrium amounts of species,
$n^{\text{q}}$; the sensitivity derivatives of this equilibrium state,
$\partial\varphi^{\text{q}}/\partial T$, $\partial\varphi^{\text{q}}/\partial P$,
and $\partial\varphi^{\text{q}}/\partial b$, needed by equation~(\ref{eq:smart-estimate});
and the derivatives of the chemical potentials, $\partial\mu^{\text{q}}/\partial T$,
$\partial\mu^{\text{q}}/\partial P$, and $\partial\mu^{\text{q}}/\partial n$,
needed by equation~(\ref{eq:mu-estimated}). See \citet{Leal2017}
for a description of how to calculate these derivatives. By saving
these state values, we will later be able to make fast and accurate
predictions of all equilibrium states in the vicinity of the one just
recorded. 

The next time the algorithm is asked to solve a new equilibrium problem,
it does so by first searching, among all saved equilibrium input conditions:
\begin{equation}
\mathcal{I}=\left\{ (T^{k},P^{k},b^{k})\right\} _{k=1}^{\text{K}},
\end{equation}
the input $(T^{\text{p}},P^{\text{p}},b^{\text{p}})$ that is closest
to the given one $(T^{\text{q}},P^{\text{q}},b^{\text{q}})$, where
``closest'' is measured here in the sense of the Euclidean norm
(other norms can be implemented as well):
\begin{equation}
d^{k}=\left[\left(T^{\text{q}}-T^{k}\right)^{2}+\left(P^{\text{q}}-P^{k}\right)^{2}+\sum_{j=1}^{\text{E}}\left(b_{j}^{\text{q}}-b_{j}^{k}\right)^{2}\right]^{\frac{1}{2}},
\end{equation}
where $d^{k}$ is a scalar that measures the difference between $(T^{\text{q}},P^{\text{q}},b^{\text{q}})$
and $(T^{k},P^{k},b^{k})$. Once the search is concluded, the previously
calculated equilibrium state, corresponding to input conditions $(T^{\text{p}},P^{\text{p}},b^{\text{p}})$,
can be used to estimate the equilibrium state corresponding to input
conditions $(T^{\text{q}},P^{\text{q}},b^{\text{q}})$ using equation~(\ref{eq:smart-estimate}).

Finally, it remains to check if the estimated equilibrium state is
potentially accurate enough using the acceptance test defined by equation~(\ref{eq:acceptance-test-mu}).
If the test succeeds, then the equilibrium calculation ends. Otherwise,
a complete chemical equilibrium calculation at $(T^{\text{q}},P^{\text{q}},b^{\text{q}})$
conditions is performed, $n^{\text{q}}=\varphi(T^{\text{q}},P^{\text{q}},b^{\text{q}})$,
and the corresponding sensitivity derivatives of the fully computed
equilibrium state ($\partial\varphi^{\text{q}}/\partial T$, $\partial\varphi^{\text{q}}/\partial P$,
and $\partial\varphi^{\text{q}}/\partial b$) are evaluated, which
are then finally saved for possible future use when estimating equilibrium
states under $(T,P,b)$ conditions in the vicinity of the just recorded
$(T^{\text{q}},P^{\text{q}},b^{\text{q}})$ conditions. 

\section{Results\label{sec:Results}}

\begin{figure}
\begin{centering}
\includegraphics[width=1\textwidth]{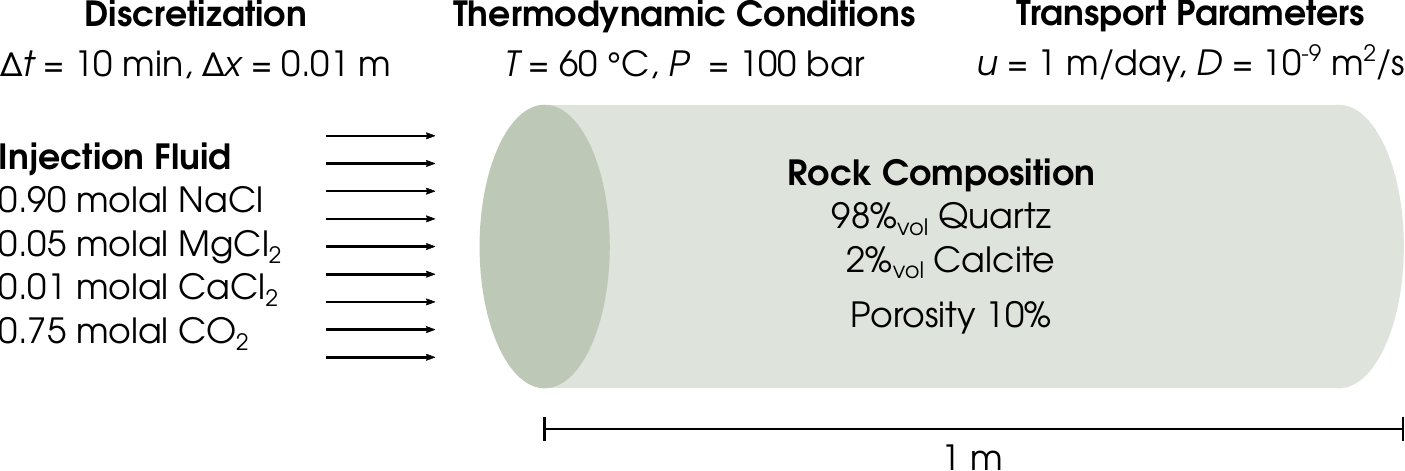}
\par\end{centering}
\caption{\label{fig:illustration-reactive-transport-model}Illustration of
the reactive transport modeling along a rock core, with details of
the injection fluid and rock composition, transport parameters, and
numerical discretization.}
\end{figure}

We now present the use of the machine learning method for smart chemical
equilibrium calculations in a reactive transport simulation and show
how its performance compares with the use of a conventional chemical
equilibrium algorithm. For details on how we solve the reactive transport
equations, see Appendix~\ref{sec:Reactive-Transport-Equations}.

Figure~\ref{fig:illustration-reactive-transport-model} illustrates
the reactive transport modeling carried out in this work to investigate
the efficiency of the proposed smart chemical equilibrium calculations
for sequential calculations. It shows the injection of an aqueous
fluid resulting from the mixture of 1 kg of water with 0.90 moles
of NaCl, 0.05 moles of MgCl$_{2}$, 0.01 moles of CaCl$_{2}$, and
0.75 moles of CO$_{2}$, in a state very close to CO$_{2}$ saturation,
and thus in an acidic state with a calculated pH of 4.65. The initial
rock composition is 98\%$_{\text{vol}}$ SiO$_{2}$(quartz) and 2\%$_{\text{vol}}$
CaCO$_{3}$(calcite), with an initial porosity of 10\%. The resident
fluid is a 0.70 molal NaCl brine in equilibrium with the rock minerals,
with a calculated pH of 10.0. The temperature and pressure of the
fluids are 60~\textdegree C and 100~bar, respectively. The reactive
transport modeling procedure assumes a constant fluid velocity of
$v=\unit[1]{m/day}$ ($\unit[1.16\cdot10^{-5}]{m/s}$) and the same
diffusion coefficient $D=\unit[10^{-9}]{m^{2}/s}$ for all fluid species,
without dispersivity.

The activity coefficients of the aqueous species are calculated using
the HKF extended Debye-Hückel model \citep{Helgeson1974,Helgeson1974a,Helgeson1976,Helgeson1981}
for solvent water and ionic species, except for the aqueous species
CO$_{2}$(aq), for which the \citet{drummond1981boiling} model is
used. The standard chemical potentials of the species are calculated
using the equations of state of \citet{Helgeson1974,Helgeson1978,Tanger1988,Shock1988}
and \citet{Shock1992}. The database file \texttt{slop98.dat} from
the software SUPCRT92 \citep{Johnson1992} is used to obtain the parameters
for the equations of state. The equation of state of \citet{Wagner2002}
is used to calculate the density of water and its temperature and
pressure derivatives. Kinetics of dissolution and precipitation of
both calcite and dolomite is neglected in this particular example
(i.e., the local equilibrium assumption is employed). 

\begin{figure}
\begin{centering}
\includegraphics[width=1\columnwidth]{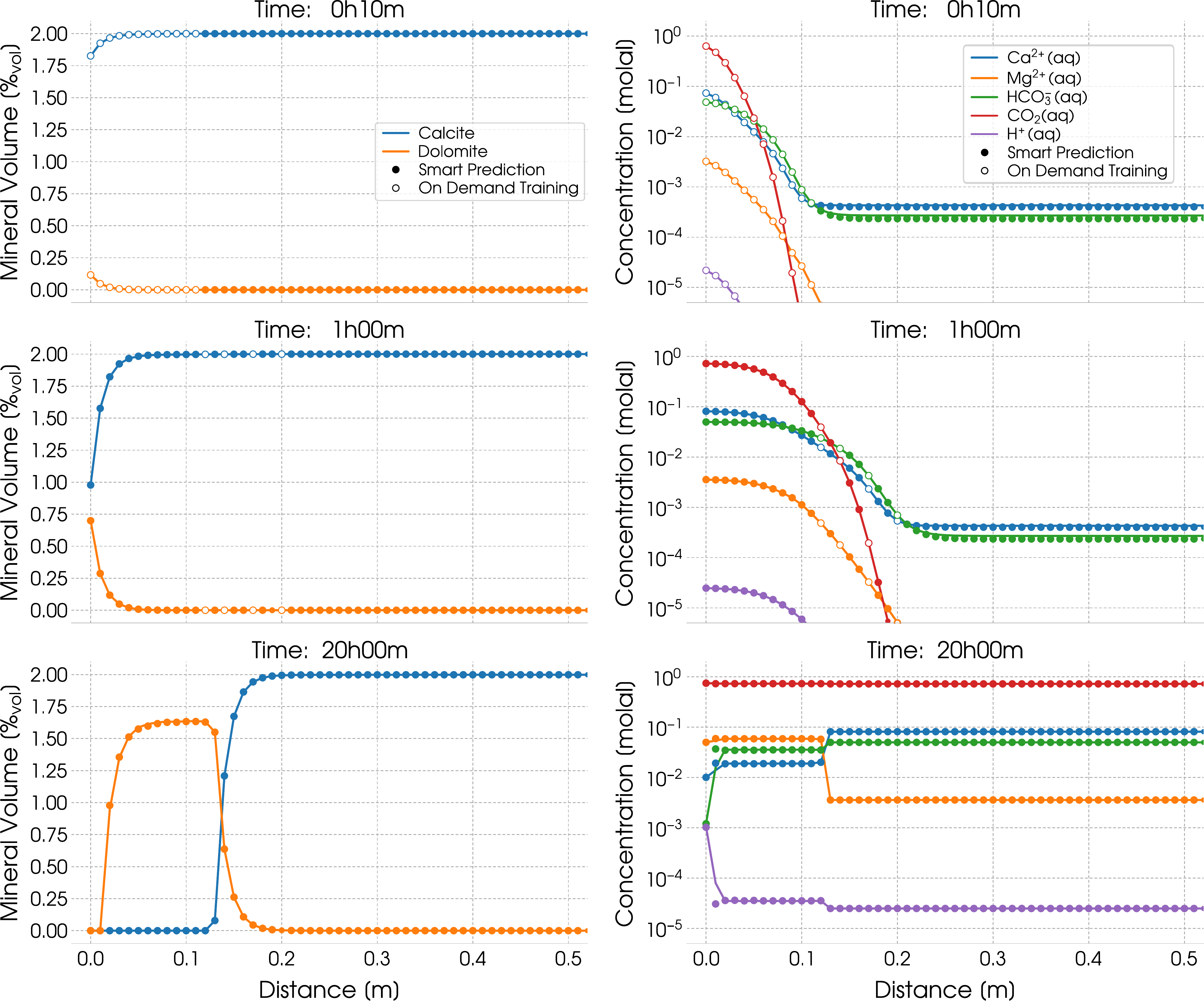}
\par\end{centering}
\caption{\label{fig:calcite-dolomite}The volume of minerals calcite and dolomite
(in $\%_{\text{vol}}$) and the concentrations of selected aqueous
species (in molal) along the rock core at three different times: 10
minutes, 1 hour, and 20 hours. The \emph{solid curves} correspond
to solving the reactive transport equations using expensive full chemical
equilibrium calculations (reference curves). The \emph{filled circles}
$\text{(}\bullet\text{)}$ denote occasions in which equilibrium states
are smartly predicted using previously recorded equilibrium states
by the machine learning algorithm, and the \emph{empty circles} $\text{(}\circ\text{)}$
denote occasions in which \emph{on-demand training} is needed by the
machine learning algorithm to learn a new chemical equilibrium state.}
\end{figure}

Figure~\ref{fig:calcite-dolomite} shows the volumes of the minerals
calcite, CaCO$_{3}$, and dolomite, CaMg(CO$_{3}$)$_{2}$, as well
as the concentrations of aqueous species Ca$^{2+}$(aq), Mg$^{2+}$(aq),
HCO$_{3}^{-}$(aq), CO$_{2}$(aq), and H$^{+}$(aq) along the rock
core at three different simulation times: 10 minutes, 1 hour, and
20 hours. As calcite dissolves, Ca$^{2+}$(aq) ions are released into
the aqueous solution, which react with the incoming Mg$^{2+}$(aq)
ions from the left boundary to precipitate dolomite. After 10 minutes
of injecting the CO$_{2}$-saturated brine, one observes a slight
dissolution of calcite and a corresponding precipitation of dolomite.
The injected CO$_{2}$-saturated brine increases the local concentrations
of carbonic species, CO$_{2}$(aq) and HCO$_{3}^{-}$(aq). The local
concentration of ions Ca$^{2+}$(aq) also increase as a result of
CaCO$_{3}$ dissolution. The precipitated dolomite, however, is gradually
dissolved, as the injection of the acidic CO$_{2}$-saturated fluid
continues. This can be seen in Figure~\ref{fig:calcite-dolomite},
in the left region of the rock core, where neither calcite nor dolomite
are present after 20 hours of continuous fluid injection. Also after
20 hours of fluid injection, the Mg$^{2+}$(aq) concentration drops
sharply between core distances 0.1~m and 0.2~m, which is exactly
where dolomite is currently precipitating. 

\begin{figure}
\begin{centering}
\includegraphics[width=1\textwidth]{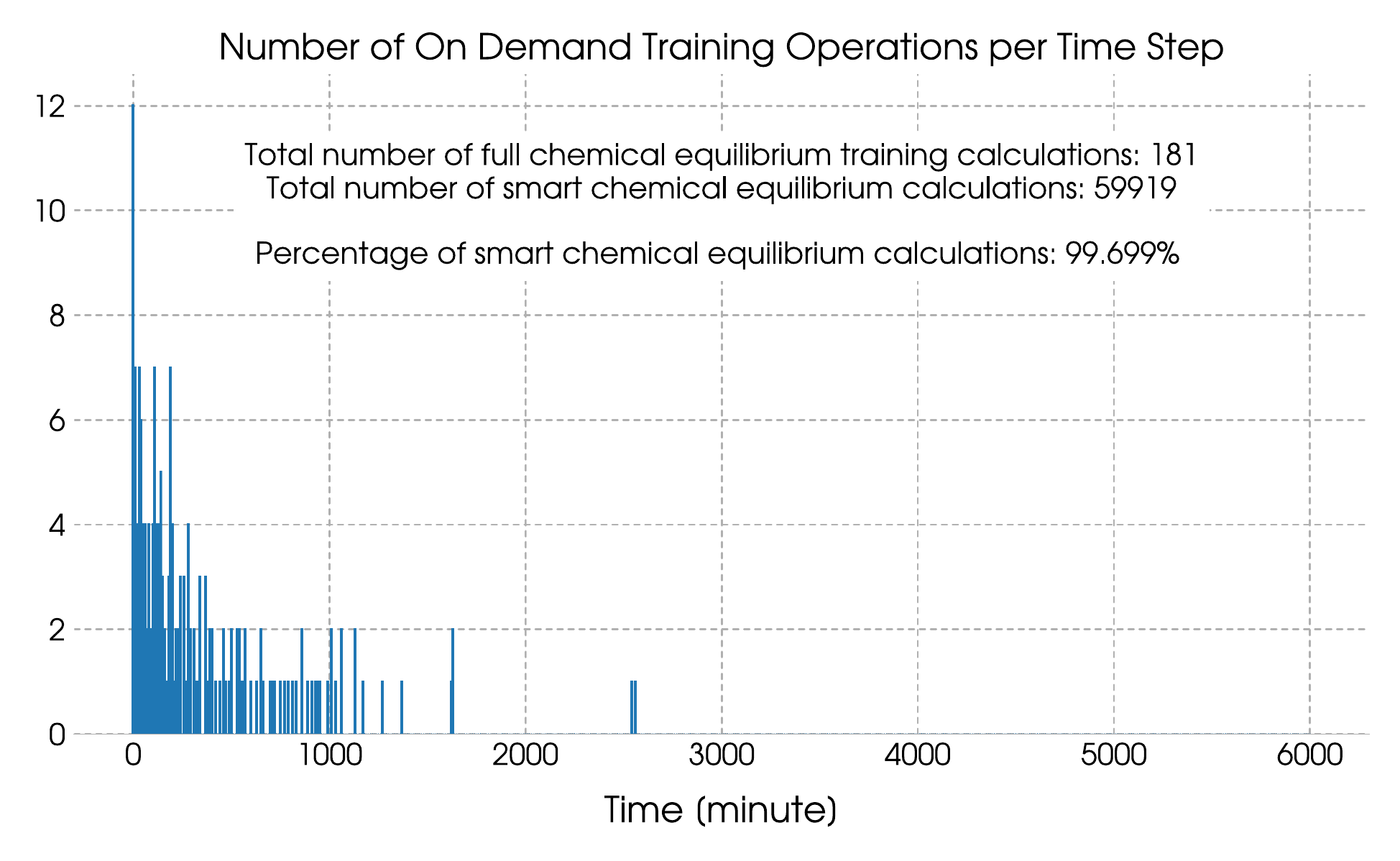}
\par\end{centering}
\caption{\label{fig:number-training}The recorded number of \emph{on-demand
training} occasions needed by the machine learning algorithm for smart
chemical equilibrium calculations during each time step (10 minutes)
of the reactive transport simulation (finished after 600 time steps,
totaling 6000 minutes of simulation). Each training required the full
solution of the non-linear equations governing chemical equilibrium
using a Newton-based numerical method \citep{Leal2016a}, which we
show here to be considerably more computationally expensive than both
transport calculations and the cost of smart equilibrium calculations
accelerated with our machine learning strategy.}
\end{figure}

Common to all sub-figures in Figure~\ref{fig:calcite-dolomite} is
the use \emph{of solid curves} to denote the reactive transport results
using conventional, and thus computationally expensive, chemical equilibrium
calculations throughout the simulation. These curves are used as a
reference case for the results obtained using the here-proposed \emph{machine
learning algorithm for smart chemical equilibrium calculations} applied
to the same reactive transport problem (see Section~\ref{subsec:Smart-chemical-equilibrium}).
Both the \emph{filled circles} $\text{(}\bullet\text{)}$ and the
\emph{empty circles} $\text{(}\circ\text{)}$ are used to mark the
calculated mineral volumes and species concentrations in Figures~\ref{fig:calcite-dolomite}
using the smart equilibrium algorithm accelerated with machine learning.
However, while the filled circles $\text{(}\bullet\text{)}$ represent
occasions, where the machine learning algorithm was able to quickly
and successfully predict an accurate equilibrium state, the empty
circles $\text{(}\circ\text{)}$ represent occasions in which the
smart algorithm required on-demand training to learn a new chemical
equilibrium problem. These on-demand learning occasions happen at
different mesh cells, either at the same time step or at different
ones. The on-demand training operation is needed, since the machine
learning algorithm is faced with an equilibrium problem for which
no accurate-enough prediction can be made. We remark that our proposed
machine learning algorithm, for smart chemical equilibrium calculations,
does not need to know any spatial or temporal information, although
this could be explored for faster search operations in the future. 

The machine learning algorithm starts without any recollection of
previous chemical equilibrium calculations. It is initially trained
during two occasions: when calculating the \emph{initial equilibrium
state} between the resident fluid and the rock minerals and when calculating
the \emph{equilibrium state of the injected fluid} at the left boundary.
We can see in Figure~\ref{fig:calcite-dolomite} that during the
first time step of the simulation, from time 0 to 10 minutes, the
smart equilibrium algorithm was able to accurately estimate the equilibrium
states in most mesh cells (see the filled circles $\bullet$). All
these successful quick estimates were done using \emph{only one out
of the two initial saved equilibrium states}: the initial equilibrium
state of fluid species and rock minerals. 

As the reactive fluid is injected inside the rock core, this promotes
strong compositional changes in both resident fluid and rock minerals.
Because of this, we can see that the machine learning algorithm required
during the first time step, at 10 minutes, \emph{on-demand training}
in the first 12 cells (see the empty circles $\circ$). As the perturbation
fronts move down the rock core, additional training is performed as
needed to fulfill a given accuracy criterion (using $\epsilon_{\text{rel}}=\epsilon_{\text{abs}}=0.1$).
As seen in Figure~\ref{fig:calcite-dolomite}, at 1h, or after 6
time steps, where strong concentration changes are occurring between
0.1~m and 0.2~m, there were 4 mesh cells that required full and
expensive chemical equilibrium calculations for both accuracy and
on-demand learning purposes. Eventually, the smart equilibrium algorithm
gains enough knowledge about the recurring equilibrium states during
the simulation, so that it is then able to quickly and accurately
perform all equilibrium calculations without further training, as
seen after 20h of fluid injection, or after 120 time steps (when only
filled circles $\bullet$ are present).

Figure~\ref{fig:number-training} shows how many on-demand training
occasions were needed at each time step. During each of these training
operations, a full chemical equilibrium calculation, using a conventional
Newton-based algorithm, as presented in \citet{Leal2016a}, is performed
and the sensitivity derivatives of the equilibrium state are calculated
for learning reasons. Figure~\ref{fig:number-training} shows some
relatively intense learning activities during the first 1000 minutes
of simulated time (or 100 time steps) as expected, since the machine
learning algorithm is constantly facing new equilibrium challenges
in the beginning of its life. During the first time step, 12 on-demand
training operations are needed, exactly in the first 12 mesh cells,
in a mesh containing 100 cells (see empty circles $\circ$ in Figure~\ref{fig:number-training}
at 10 minutes). This implies that during the first time step, the
machine learning algorithm was not yet trained enough to accurately
predict the equilibrium states in those cells. During the second time
step, 7 additional on-demand training operations are carried out and,
as the simulation continues, only 1–2 additional training cases happen
subsequently per time step, until the smart chemical equilibrium algorithm
no longer requires any further learning (after about 258 time steps
of 10 minutes of simulated time each). At this point, the machine
learning algorithm is able to quickly and accurately predict all subsequent
equilibrium states in every mesh cell, at every time step. For the
reactive transport modeling problem chosen here, from the total of
60,000 equilibrium calculations needed, only 181 of these were actually
performed using the expensive Newton-based chemical equilibrium algorithm,
so that 99.7\% of all such calculations were quickly performed using
the proposed machine learning algorithm for fast chemical equilibrium
calculations.

\begin{figure}
\begin{centering}
\includegraphics[width=1\textwidth]{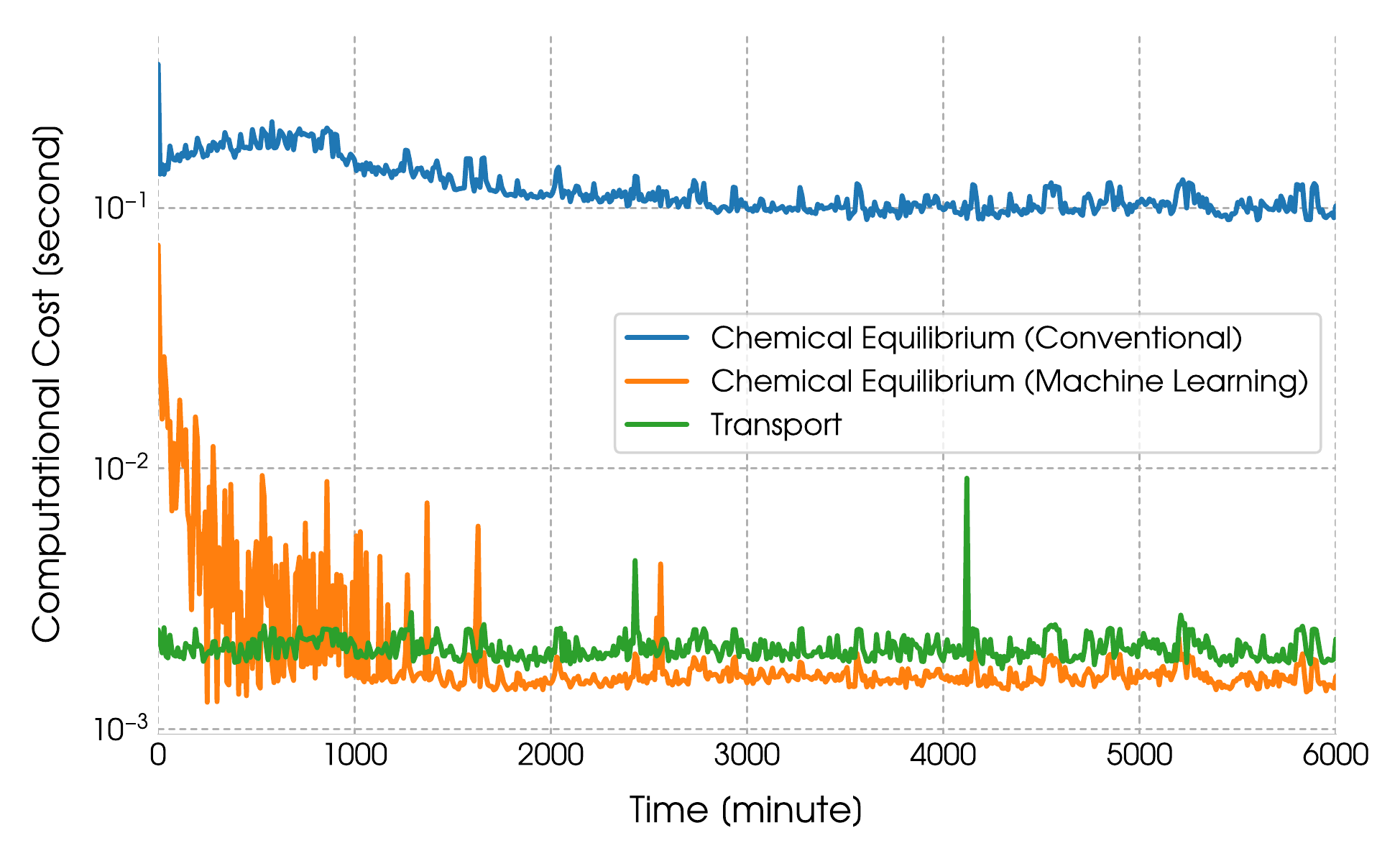}
\par\end{centering}
\caption{\label{fig:computational-cost}The comparison of the computational
cost of transport, conventional chemical equilibrium, and machine-learning-accelerated
chemical equilibrium calculations in each time step of the reactive
transport simulation. The computational cost is given as CPU time
(in seconds). The cost of equilibrium calculations per time step is
the sum of the individual cost in each mesh cell. The cost of transport
calculations per time step is the cost of solving the discretized
algebraic transport equations. Our proposed machine learning strategy,
based on predictions using sensitivity derivatives, can significantly
speedup chemical equilibrium calculations by almost two orders of
magnitude in this specific example.}
\end{figure}

Figure~\ref{fig:computational-cost} compares the computational cost
per time step for conventional and smart chemical equilibrium calculations
with the cost for transport calculations. The computational costs
are measured in CPU time (in seconds). For the equilibrium calculations,
it is the CPU time taken to calculate the equilibrium states of all
mesh cells in a time step. For the transport calculations, it is the
CPU time taken to solve the algebraic transport equations (see Appendix~\ref{sec:Reactive-Transport-Equations}).
Figure~\ref{fig:computational-cost} shows that the cost for full,
conventional chemical equilibrium calculations is about two orders
of magnitude greater than the cost for transport calculations. It
also shows that the computational cost for the machine learning algorithm,
for smart equilibrium calculations, is greater than that of transport
calculations \emph{at the beginning}, but less than that for conventional
Newton-based equilibrium calculations. After the initial on-demand
training phase of the machine-learning-accelerated chemical equilibrium
algorithm, it can be seen in Figure~\ref{fig:computational-cost}
that its computational cost stabilizes and remains below the cost
of transport calculations. To note, typically, full chemical equilibrium
calculations require by far the most computing time during reactive
transport simulations, as also seen in Figure~\ref{fig:computational-cost}
and as stated in the introduction.

\begin{figure}
\begin{centering}
\includegraphics[width=1\textwidth]{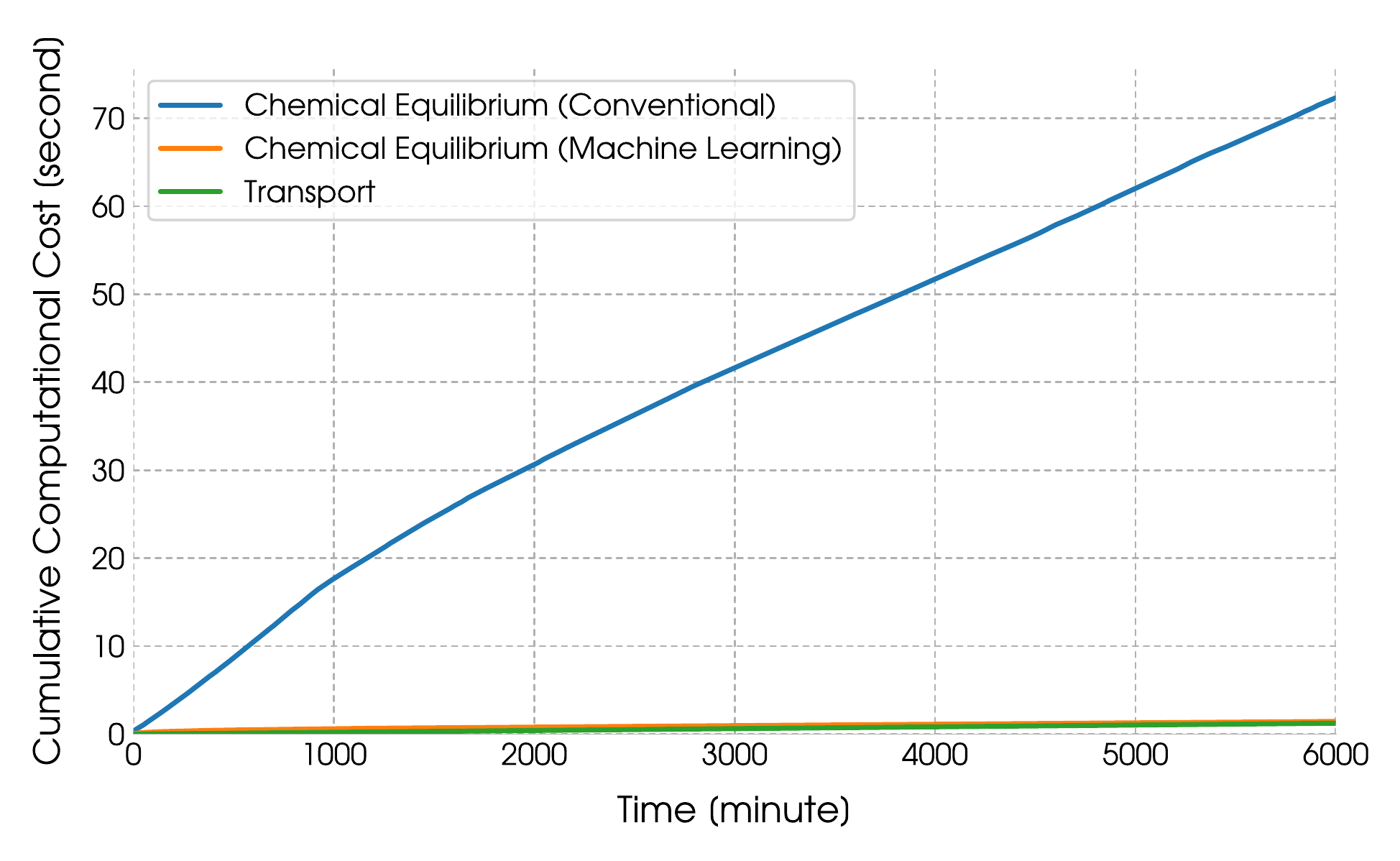}
\par\end{centering}
\caption{\label{fig:computational-cost-cumulative}The cumulative computational
cost of transport, conventional chemical equilibrium, and machine-learning-accelerated
chemical equilibrium calculations as the reactive transport simulation
continues. Our proposed smart chemical equilibrium algorithm promotes
substantial savings in computational costs, resulting in comparable
costs to that of transport.}
\end{figure}

Figure~\ref{fig:computational-cost-cumulative} compares the \emph{cumulative
computational costs} of conventional and machine-learning-accelerated
chemical equilibrium calculations with the cost of transport calculations.
After several time steps (about 200 or at about 2000 minutes of simulated
time in the figure), the cumulative computational costs for all these
calculations increase at a more or less constant rate. However, it
can be seen that the increment rate of the cumulative costs for the
conventional equilibrium calculations, using a Newton-based method,
is considerably higher than that for transport (about 45 times higher)
and than that for smart, machine-learning-accelerated equilibrium
calculations (about 63 times higher).

\begin{figure}
\begin{centering}
\includegraphics[width=1\textwidth]{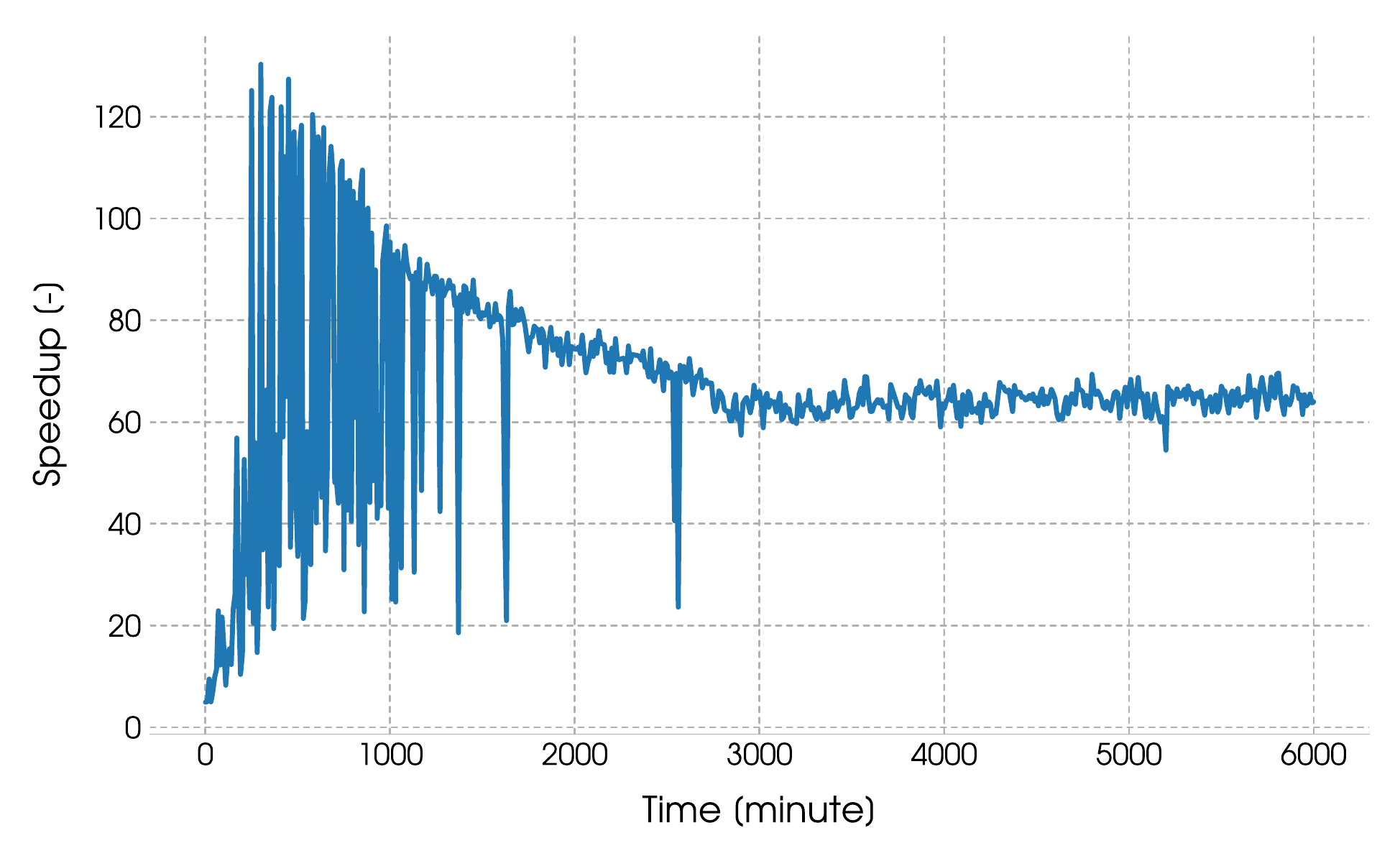}
\par\end{centering}
\caption{\label{fig:smar-vs-full-algorithm-speedup}The speedup of chemical
equilibrium calculations, at each time step of the simulation, when
employing our proposed machine learning acceleration strategy for
fast and accurate predictions of chemical equilibrium states. The
speedup, in a time step, is calculated as the ratio of the time needed
for the conventional chemical equilibrium algorithm and the time needed
for the smart chemical equilibrium algorithm. \textbf{Remark:} By
implementing more optimized lookup algorithms, when finding the most
similar previously solved equilibrium problem (e.g., using kd-trees,
which permits a search algorithm with complexity of $O(\log n)$,
instead of the naive search of complexity $O(n)$ used here, where
$n$ is the number of saved equilibrium inputs), the speedup shown
here will likely increase significantly. }
\end{figure}

Figure~\ref{fig:smar-vs-full-algorithm-speedup} shows the speedup
during chemical equilibrium calculations, at each time step of the
simulation, promoted by the use of the smart chemical equilibrium
algorithm. The speedup at a given time step is calculated as the ratio
of the CPU time needed for conventional equilibrium calculations and
the CPU time for smart equilibrium calculations, performed in all
mesh cells at that time step. At the beginning, during a more active
on-demand training phase by the smart chemical equilibrium solver,
the speedup oscillates, achieving a maximum of about 125. The oscillation
is a combination of many factors, such as the need to occasionally
perform full equilibrium calculations for learning purposes and the
fact that the strongest compositional changes occur during the first
time steps when the injected fluid perturbs its initial equilibrium.
As the machine learning chemical equilibrium algorithm continues to
learn how to solve different equilibrium problems, it gets to a stage
that is able to quickly predict all new equilibrium states. When this
happens, the speedup stabilizes at about 63 (or $10^{1.8}$), as shown
in Figure~\ref{fig:smar-vs-full-algorithm-speedup}, at later time
steps. 

We remark that this speedup is strongly dependent on how fast we can
search for the previous equilibrium problem, closest to the one being
solved. Currently, we are using a naive approach, during which all
past equilibrium inputs, $\mathcal{I}=\{(T^{k},P^{k},b^{k})\}_{k=1}^{\text{K}}$,
are saved in a linked list data structure. This data structure is
not optimal for \emph{nearest-neighbor search algorithms}, a common
procedure in computer science, such as in machine learning and computer
graphics applications, because it requires each previous input $(T^{k},P^{k},b^{k})$
to be compared against the new equilibrium input $(T^{\text{q}},P^{\text{q}},b^{\text{q}})$.
As a result, this particular search algorithm has complexity $O(\text{K})$,
where K is the number of saved inputs. A more suitable data structure
is a \emph{kd-tree}, which would split the input space into several
orthogonal partitions, permitting a search algorithm that succeeds,
on average, after about $\log_{2}(\text{K})$ comparisons. For example,
instead of 256 comparisons using a linked list, we would need only
$8=\log_{2}(256)$ with a balanced kd-tree. Future optimization of
our machine learning equilibrium method will make use of kd-trees
or other fast search algorithms. 

\section{Conclusions and Future Work\label{sec:Conclusions-and-Future-Work}}

We present a straightforward, smart, machine-learning chemical equilibrium
algorithm that calculates equilibrium states in reactive transport
simulations 60–125 times faster than a conventional Newton-based
algorithm. The reason for this speedup is that the smart equilibrium
algorithm employs a new concept of an unconventional, supervised machine
learning method, with an on-demand training strategy, rather than
the more common training\nobreakdash-in\nobreakdash-advance approach.
Our new algorithm is capable of ``remembering'' past equilibrium
calculations and performing rapid predictions of new equilibrium states,
whenever the new equilibrium problem is sufficiently close to some
previously solved one. These machine-learning-accelerated chemical
equilibrium calculations are able to dramatically speedup reactive
transport simulations, in which chemical reaction calculations are
otherwise typically responsible for over 90\% of all computation costs.

Our proposed machine learning algorithm differs from previous machine
learning strategies to speedup chemical equilibrium calculations,
because here, learning is carried out on-demand, i.e., during the
actual simulation, rather than during an initial, time-consuming training
phase. Consequently, our method requires no a priori insights of possible
chemical conditions that may occur during the actual simulation, a
problem inherent to standard, i.e., statistics-based, machine learning
algorithms. These initial training sessions of standard machine learning
algorithms require large data-sets of input-output information, which
are sometimes randomly obtained. Furthermore, there are no clear criteria
for determining when such training sessions are sufficient and can
be terminated to begin with the actual simulation of interest. Finally,
traditional, statistically\nobreakdash-based machine learning algorithms
applied to chemical equilibrium calculations neither understand the
thermodynamic behavior of stable phases in chemical systems nor can
they straightforwardly predict chemical equilibrium states that satisfy
given mass conservation constraints. 

The presented smart equilibrium algorithm contains only three steps,
which are computationally relatively cheap to perform. The first step
performs a lookup of past equilibrium problems that are closest to
the new one. The second step is a simple algebraic calculation to
estimate the new equilibrium state from the previous one, more specifically,
two vector additions (of dimension $\text{N}\times1$) and a matrix-vector
multiplication (the matrix with dimension $\text{N}\times\text{E}$,
the vector with dimension $\text{E\ensuremath{\times}1}$), where
N and E are the number of species and elements, respectively. The
third step is to check the acceptance criterion, which is also a relatively
cheap compute operation. Hence, the smart approximation is able to
quickly bypass extremely expensive operations, when solving chemical
equilibrium problems, such as \emph{(i)} the evaluation of thermodynamic
properties of all species (e.g., activities)and \emph{(ii)} the solution
of systems of linear equations. 

Since equilibrium calculations are carried out iteratively, these
computationally expensive operations are performed once each iteration
and, thus, several times. Even if convergence could always be established
in only one iteration, this would still not cause a substantial improvement
in the performance of reactive transport simulations, given the high
cost of those operations. In contrast, the proposed smart equilibrium
algorithm \emph{skips all these iterations,} and their intrinsic expensive
operations, and is instead able to immediately predict an accurate
equilibrium state for similar problems. The new algorithm, thus, has
the potential to substantially accelerate reactive transport simulations,
as, indeed, demonstrated for the modeling problem in Section~\ref{sec:Results}.

Potential future work and improvements of our machine learning algorithm
for rapid chemical equilibrium calculations could include:
\begin{itemize}
\item the use of kd-trees for faster lookup operations;
\item the investigation of the smart chemical equilibrium algorithm for
more complex chemical systems in more heterogeneous and complex reactive
transport problems per time step;
\item the use of alternative acceptance criteria to further avoid unnecessary
on-demand training operations;
\item the use of GPUs for massively parallel smart predictions of thousands
to millions of equilibrium states;
\item the extension of the machine learning strategy, presented here, to
chemical kinetics calculations and other problems in science and engineering;
\item the application of \emph{garbage collector ideas}, used in programming
languages, such as Python and Java, to eliminate all previously saved
equilibrium states that have not been used for some time;
\item the implementation of strategies that rely not only on the nearest
previous saved equilibrium state, but up to a certain number (e.g.,
the two or three closest equilibrium problems solved previously) to
possibly avoid the triggering of on-demand training operations and
also to improve accuracy by combining the sensitivity derivatives
in different nearby states;
\item incorporate spatial and temporal information in the machine learning
algorithm and test if the computation cost of search operations can
be decreased with it.
\end{itemize}
The above roadmap of future investigations shows how much a new machine
learning algorithm for rapid chemical reaction calculations can yet
be improved in addition to what we have already shown here. We believe
that the use of machine\nobreakdash-learning\nobreakdash-accelerated
algorithms for chemical reaction calculations is crucial for a significant
acceleration of complex reactive transport simulations, and it is
essential for a substantial decrease of the overall computational
costs of chemical calculations, which so far have been responsible
for 90–99\% of all computing costs in large-scale numerical simulations
with intricate chemistry representation. 

\bibliographystyle{apalike-order-by-citation}
\bibliography{library}

\appendix

\section{Chemical Equilibrium Equations \label{subsec:Chemical-equilibrium-equations}}

The solution of the Gibbs energy minimization problem in equation~(\ref{eq:gem-problem})
needs to satisfy the following \emph{Karush–Kuhn–Tucker (KKT) conditions},
or \emph{first-order optimality conditions,} for a local minimum of
the Gibbs energy function $G$ \citep{Nocedal1999,Fletcher2000}:
\begin{alignat}{2}
\mu-A^{T}y-z & =0,\\
An-b & =0,\\
n_{i}z_{i} & =0 & \qquad & (i=1,\ldots,\text{N}),\\
n_{i} & \geq0 &  & (i=1,\ldots,\text{N}),\\
z_{i} & \geq0 &  & (i=1,\ldots,\text{N}),
\end{alignat}
where ${y=(y_{1},\ldots,y_{\text{E}})}$ and ${z=(z_{1},\ldots,z_{\text{N}})}$
are introduced \emph{Lagrange multipliers} that need to be solved
along with the species amounts ${n=(n_{1},\ldots,n_{\text{N}})}$.
For more details about these Lagrange multipliers and their interpretation
as well as instructions on how to efficiently solve these equations,
see \citet{Leal2016a,Leal2017}. 

The previous chemical equilibrium equations can be written in an \emph{extended
law of mass action} (xLMA) formulation as:
\begin{alignat}{2}
\ln K-\nu(\ln a+\ln w) & =0,\\
An-b & =0,\\
n_{i}\ln w_{i} & =0 & \qquad & (i=1,\ldots,\text{N}),\\
n_{i} & \geq0 &  & (i=1,\ldots,\text{N}),\\
0<w_{i} & \leq1 &  & (i=1,\ldots,\text{N})
\end{alignat}
following the use of the extended law of mass action equations:
\begin{equation}
K_{m}=\prod_{i=1}^{\text{N}}(a_{i}w_{i})^{\nu_{mi}}\label{eq:xlma-equations}
\end{equation}
associated with the M linearly independent chemical reactions among
the N chemical species in equilibrium: 
\begin{equation}
0\rightleftharpoons\sum_{i=1}^{\text{N}}\nu_{mi}\alpha_{i}\qquad(m=1,\ldots,\text{M}),
\end{equation}
where ${K=(K_{1},\ldots,K_{\text{M}})}$ is the vector of \emph{equilibrium
constants} of the reactions, with ${K_{m}=K_{m}(T,P)}$ denoting the
equilibrium constant of the $m$th reaction, and $\nu$ the $\text{M}\times\text{N}$
\emph{stoichiometric matrix} of these chemical reactions, with $\nu_{mi}$
denoting the stoichiometric coefficient of the $i$th species in the
$m$th reaction with the convention that $\nu_{mi}$ is positive if
the $i$th species is a product in the $m$th reaction, and negative
if a reactant. Moreover, ${w=(w_{1},\ldots,w_{\text{N}})}$ is the
vector of \emph{species stability factors} that need to be solved
along with the species amounts ${n=(n_{1},\ldots,n_{\text{N}})}$.
These factors are introduced to ensure that the extended law of mass
action equations (\ref{eq:xlma-equations}) are valid even when some
species in their corresponding reactions are unstable at equilibrium
(i.e., when a species belongs to a phase that is absent from equilibrium).
When all species are stable at equilibrium, it follows that ${w_{i}=1}$
and the xLMA equations reduce to the conventional LMA equations:
\begin{equation}
K_{m}=\prod_{i=1}^{\text{N}}a_{i}^{\nu_{mi}}.
\end{equation}
The number of linearly independent chemical reactions among the N
species in equilibrium is ${\text{M}=\text{N}-\text{C}}$, where ${\text{C}=\text{rank}(A)}$,
and thus whenever the formula matrix $A$ is full rank, $\text{C}=\text{E}$
and $\text{M}=\text{N}-\text{E}$. For more information on how to
solve these equations and how they are related to the conventional
law of mass action equations, see \citet{Leal2016c,Leal2017}.

\section{Reactive Transport Equations\label{sec:Reactive-Transport-Equations}}

The fundamental mass conservation equations for both fluid and solid
species are:
\begin{alignat}{2}
\frac{\partial n_{i}^{\text{f}}}{\partial t}+\nabla\cdot(\boldsymbol{v}n_{i}^{\text{f}}-D\nabla n_{i}^{\text{f}}) & =r_{i}^{\text{f}} & \qquad & (i=1,\ldots,\text{N}^{\text{f}}),\label{eq:cons-mass-fluid-species}\\
\frac{\partial n_{i}^{\text{s}}}{\partial t} & =r_{i}^{\text{s}} &  & (i=1,\ldots,\text{N}^{\text{s}}),\label{eq:cons-mass-solid-species}
\end{alignat}
where $n_{i}^{\text{f}}$ and $n_{i}^{\text{s}}$ are the \emph{bulk
concentration} of the $i$th fluid and solid species (in mol/m$^{3}$),
respectively; $\boldsymbol{v}$ is the fluid pore velocity (in m/s);
$D$ is the diffusion coefficient of the fluid species (in m$^{2}$/s);
$r_{i}^{\text{f}}$ and $r_{i}^{\text{s}}$ are the rates of production\slash{}consumption
of the $i$th fluid and solid species (in mol\slash{}s), respectively,
due to chemical reactions; and $\text{N}^{\text{f}}$ and $\text{N}^{\text{s}}$
are the number of fluid and solid species, respectively. Note that
the above equations assume a single fluid phase and common diffusion
coefficients for all fluid species.

By partitioning the species as fluid and solid species, the formula
matrix, $A,$ can be conveniently represented as:
\begin{equation}
A=\begin{bmatrix}A^{\text{f}} & A^{\text{s}}\end{bmatrix},
\end{equation}
where $A^{\text{f}}$ and $A^{\text{s}}$ are the formula matrices
of the fluid and solid partitions (i.e., the matrices constructed
from the columns of $A$ corresponding to fluid and solid species).
The concentrations of elements in both fluid and solid partitions,
$b_{j}^{\text{f}}$ and $b_{j}^{\text{s}}$, can then be calculated
from the species concentrations in the same partition using: 
\begin{alignat}{2}
b_{j}^{\text{f}} & =\sum_{i=1}^{\text{N}^{\text{f}}}A_{ji}^{\text{f}}n_{i}^{\text{f}} & \qquad & (j=1,\ldots,\text{E})\\
\shortintertext{and}b_{j}^{\text{s}} & =\sum_{i=1}^{\text{N}^{\text{s}}}A_{ji}^{\text{s}}n_{i}^{\text{s}} &  & (j=1,\ldots,\text{E}).
\end{alignat}

Recall that the rates of production of the species, $r_{i}^{\text{f}}$
and $r_{i}^{\text{s}}$, are exclusively due to chemical reactions.
Let $r_{i}$ denote the rate of production\slash{}consumption of
the $i$th species in the system, and not in a fluid\slash{}solid
partition. From the mass conservation condition for the elements (chemical
elements and electrical charge), it follows that:
\begin{equation}
\underset{\substack{\text{rate of production}\\
\text{ of element \ensuremath{j}}
}
}{\underbrace{\sum_{i=1}^{\text{N}}A_{ji}r_{i}}}=\underset{\substack{\text{rate of production}\\
\text{ of element \ensuremath{j} in}\\
\text{the fluid partition}
}
}{\underbrace{\sum_{i=1}^{\text{N}^{\text{f}}}A_{ji}^{\text{f}}r_{i}^{\text{f}}}}+\underset{\substack{\text{rate of production}\\
\text{ of element \ensuremath{j} in}\\
\text{the solid partition}
}
}{\underbrace{\sum_{i=1}^{\text{N}^{\text{s}}}A_{ji}^{\text{s}}r_{i}^{\text{s}}}}=\quad0,
\end{equation}
which is the mathematical statement for the fact that \emph{elements
are neither created nor destroyed} during chemical reactions. We can
combine this result with equations (\ref{eq:cons-mass-fluid-species})
and (\ref{eq:cons-mass-solid-species}) to derive the following conservation
equations for the elements:
\begin{equation}
\frac{\partial b_{j}^{\text{s}}}{\partial t}+\frac{\partial b_{j}^{\text{f}}}{\partial t}+\nabla\cdot(\boldsymbol{v}b_{j}^{\text{f}}-D\nabla b_{j}^{\text{f}})=0\qquad(j=1,\ldots,\text{E}).
\end{equation}

Assume that all species, fluid and solid, are in \emph{local chemical
equilibrium everywhere, at all times}. One can then perform \emph{operator
splitting steps} to solve the fundamental mass conservation equations
(\ref{eq:cons-mass-fluid-species}) and (\ref{eq:cons-mass-solid-species})
to calculate the concentrations of the species, $n_{i}$, over time.
Let $k$ denote the current \emph{time step} and $\Delta t$ the \emph{time
step length} used in the discretization of the time derivative terms.
The operator splitting steps at the $k$th time step are:
\begin{enumerate}[wide=0.0\parindent,label=\bf{Step \arabic*)}]
\item update the concentrations of the elements in the fluid partition
using:
\[
\frac{\tilde{b}_{j}^{\text{f},k+1}-b_{j}^{\text{f},k}}{\Delta t}+\nabla\cdot(\boldsymbol{v}\tilde{b}_{j}^{\text{f},k+1}-D\nabla\tilde{b}_{j}^{\text{f},k+1})=0\qquad(j=1,\ldots,\text{E}),
\]
where $\tilde{b}_{j}^{\text{f},k}$ is the known concentration of
the $j$th element at time step $k$ in the fluid partition and $\tilde{b}_{j}^{\text{f},k+1}$
is its unknown concentration at the next time step. Note that an implicit
scheme is assumed for both advection and diffusion rates.
\item update the total concentrations of the elements:
\begin{equation}
b_{j}^{k+1}=\tilde{b}_{j}^{\text{f},k+1}+b_{j}^{\text{s},k}.
\end{equation}
\item calculate the concentrations of the species, $n_{i}^{k+1}$, in each
mesh cell, using a conventional or the proposed smart chemical equilibrium
algorithm. For this, use the local temperature and pressure values
together with the \emph{updated local concentrations of elements},
$b_{j}^{k+1}$ as inputs.
\end{enumerate}

\end{document}